%common header stuff

\documentstyle{amsppt}
\baselineskip18pt
\magnification=\magstep1
%\NoPageNumbers
%\NoRunningHeads
%\pagewidth{4.5in}
%\pageheight{7.0in}
\pagewidth{30pc}
\pageheight{45pc}
\hyphenation{co-deter-min-ant co-deter-min-ants pa-ra-met-rised
pre-print pro-pa-gat-ing pro-pa-gate
fel-low-ship Cox-et-er dis-trib-ut-ive}
\def\leaderfill{\leaders\hbox to 1em{\hss.\hss}\hfill}

\def\C{{\Cal C}}

\

\def\idest{i.e.,\ }

\def\a{{\alpha}}
\def\be{{\beta}}
\def\g{{\gamma}}

\def\e{{\varepsilon}}

\def\l{{\lambda}}

\def\b0{\text{\bf 0}}

\def\ra{{\ \longrightarrow \ }}
\def\sra{{\rightarrow}}

\def\zed{{\Bbb Z}}

\def\boxit#1{\vbox{\hrule\hbox{\vrule \kern3pt
\vbox{\kern3pt\hbox{#1}\kern3pt}\kern3pt\vrule}\hrule}}
\def\rabbit{\vbox{\hbox{\kern0pt
\vbox{\kern0pt{\hbox{---}}\kern3.5pt}}}}

\def\tableau#1{
        \hbox {
                \hskip -10pt plus0pt minus0pt
                \raise\baselineskip\hbox{
                \offinterlineskip
                \hbox{#1}}
                \hskip0.25em
        }
}

\def\tabCol#1{
\hbox{\vtop{\hrule
\halign{\strut\vrule\hskip0.5em##\hskip0.5em\hfill\vrule\cr\lower0pt
\hbox\bgroup$#1$\egroup \cr}
\hrule
} } \hskip -10.5pt plus0pt minus0pt}

\def\CR{
        $\egroup\cr
        \noalign{\hrule}
        \lower0pt\hbox\bgroup$
}

% Set up the map arrows for commutative diagrams.
%\def\mapright#1{\smash{
%     \mathop{\longrightarrow}\limits^{#1}}}

%Set up macro for commutative diagrams etc. (see Ex. 18.46 in TeXbook)

\def\blank#1#2{
%\hbox to {{#1}}{\vbox to {{#2}}}
\hbox to #1{\hfill \vbox to #2{\vfill}}
}

%Ross's table macros

\def\strut{\vrule height10pt depth5pt width0pt}

\topmatter
\title On rank functions for heaps
\endtitle

\author R.M. Green \endauthor
\affil 
Department of Mathematics and Statistics\\ Lancaster University\\
Lancaster LA1 4YF\\ England\\
{\it  E-mail:} r.m.green\@lancaster.ac.uk\\
\endaffil

\abstract
Motivated by work of Stembridge, we study rank functions for Viennot's 
heaps of pieces.  We produce a simple and sufficient criterion for a heap to
be a ranked poset and apply the results to the heaps arising from fully 
commutative words in Coxeter groups.
\endabstract


\endtopmatter

\centerline{\bf To appear in the Journal of Combinatorial Theory,
Series A}

%fixed error in Lemma 1.2.8 following second referee report: 24/2/03.
%revised in light of referee comments: 13/2/03.  Also changed 'closed
%interval' to 'balanced interval'.
%JCTA submittable version.  26/6/02.

\head Introduction \endhead

A heap is an isomorphism class of labelled posets satisfying certain 
axioms.  Heaps have a wide variety of applications, as discussed by Viennot 
in \cite{{\bf 7}}.  Stembridge \cite{{\bf 5}} showed how to associate 
heaps to fully commutative elements of Coxeter groups; the latter are the 
elements for which any
reduced expression may be obtained from any other by iterated commutation of
adjacent Coxeter generators.  In \cite{{\bf 6}}, Stembridge applied 
these ideas to $\l$-minuscule elements of Coxeter groups; these were 
first introduced by D. Peterson (unpublished) and were shown to be fully 
commutative by Proctor \cite{{\bf 4}}.

It follows from \cite{{\bf 6}, Corollary 3.4} that, under the extra assumption
that the labels occurring in the heap index an acyclic subset of the Coxeter
graph, the heap of a minuscule element is ranked as an abstract poset.
In the light of this result, it is natural to ask under what circumstances
a heap is ranked, and furthermore, what can be said about the case of heaps of
fully commutative elements of Coxeter groups?  We maintain the assumption of
\cite{{\bf 6}, Corollary 3.4}---because, as we explain in \S2.1, the situation 
becomes much more complicated otherwise---and we obtain in Theorem 2.1.1 
a simple necessary and sufficient condition for a heap to be ranked, which
involves the consideration of certain subintervals.  
We also look in \S3 at the special case of heaps of fully commutative 
elements of finite Coxeter groups, where our necessary and sufficient 
condition can be refined so that it is more explicit and easier to 
verify (Theorem 3.2.3).  For a Coxeter group of type $A$, the situation 
is simpler still and our main results are already known in this case 
(see Remark 3.3.7).

In the computer science literature \cite{{\bf 3}}, 
heaps have been used to model concurrency,
where the elements of the heap represent processes.  It would be interesting
to know if rank functions for heaps have implications for the scheduling of
such processes.

\head 1. Preliminaries \endhead

\subhead 1.1 Heaps \endsubhead

We start by recalling the basic definitions associated to heaps.  Our notation
largely follows that of \cite{{\bf 7}}.

\definition{Definition 1.1.1}
Let $P$ be a set equipped with a symmetric and reflexive binary relation
$\C$.  The elements of $P$ are called {\it pieces}, and the relation
$\C$ is called the {\it concurrency relation}.

A {\it labelled heap} with pieces in $P$ is a triple $(E, \leq, \e)$ 
where $(E, \leq)$ is a finite (possibly empty)
partially ordered set with order relation denoted
by $\leq$ and $\e$ is a map $\e : E \ra P$ satisfying the following two
axioms. 

\item{1.}{For every $\a, \be \in E$ such that $\e(\a) \ \C \ \e(\be)$, 
$\a$ and $\be$ are comparable in the order $\leq$.}

\item{2.}{The order relation $\leq$ is the transitive closure of the
relation $\leq_\C$ such that for all $\a, \be \in E$, $\a \ \leq_\C \ \be$ 
if and only if both $\a \leq \be$ and $\e(\a) \ \C \ \e(\be)$.}
\enddefinition

The terms {\it minimal} and {\it maximal} applied to the elements of 
the labelled heap refer to minimality (respectively, maximality) with 
respect to $\leq$.

\example{Example 1.1.2}
Let $P = \{1, 2, 3\}$ and, for $x, y \in P$, define $a \ \C \ b$ 
if and only if 
$|x - y| \leq 1$.  Let $E = \{a, b, c, d, e\}$ partially ordered by
extension of the relations $a \leq c$, $b \leq c$, $c \leq d$,
$c \leq e$.  Define the map $\e$ by the conditions $\e(a) = \e(d) = 1$,
$\e(c) = 2$ and $\e(b) = \e(e) = 3$.  Then $(E, \leq, \e)$ can easily be
checked to satisfy the axioms of Definition 1.1.1 and it is a labelled heap.
The minimal elements are $a$ and $b$, and the maximal elements are $d$ and $e$.
\endexample

\definition{Definition 1.1.3}
Let $(E, \leq, \e)$ and $(E', \leq', \e')$ be two labelled 
heaps with pieces in $P$ and with the same concurrency relation, $\C$.  
Two labelled heaps are isomorphic if there is a poset isomorphism
$\phi : E \ra E'$ such that $\e = \e' \circ \phi$ (\idest a labelled
poset isomorphism).

A {\it heap} of pieces in $P$ with concurrency relation $\C$ is a labelled
heap (Definition 1.1.1) defined up to labelled poset isomorphism.
The set of such heaps is denoted by $H(P, \C)$.  We denote the heap 
corresponding to the labelled heap $(E, \leq, \e)$ by $[E, \leq, \e]$.
\enddefinition

We will sometimes abuse language and speak of the underlying set of a heap,
when what is meant is the underlying set of one of its representatives.

\definition{Definition 1.1.4}
Let $(E, \leq, \e)$ be a labelled heap with pieces in $P$ and $F$ a 
subset of $E$.  
Let $\e'$ be the restriction of $\e$ to $F$.  Let ${\Cal R}$ be the relation
defined on $F$ by $\a \ {\Cal R} \ \be$ if and only if $\a \leq \be$ and
$\e(\a) \ \C \ \e(\be)$.  Let $\leq'$ be the transitive closure of ${\Cal R}$.
Then $(F, \leq', \e')$ is a labelled heap with pieces in $P$.  The heap 
$[F, \leq', \e']$ is called a {\it subheap} of $[E, \leq, \e]$.
\enddefinition

We will often implicitly use the fact that a subheap is determined by its
set of vertices and the heap it comes from.

\definition{Definition 1.1.5}
Let $E = [E, \leq_E, \e]$ and $F = [F, \leq_F, \e']$ be two heaps in 
$H(P, \C)$.
We define the heap $G = [G, \leq_G, \e''] = E \circ F$ of $H(P, \C)$ 
as follows.

\item{1.}{The underlying set $G$ is the disjoint union of $E$ and $F$.}
\item{2.}{The labelling map $\e''$ is the unique map $\e'' : G \ra P$ whose
restriction to $E$ (respectively, $F$) is $\e$ (respectively, $\e'$).}
\item{3.}{The order relation $\leq_G$ is the transitive closure of the
relation ${\Cal R}$ on $G$, where $\a \ {\Cal R} \ \be$ if
and only if one of the following three conditions holds:
\item{(i)}{$\a, \be \in E$ and $\a \leq_E \be$;}
\item{(ii)}{$\a, \be \in F$ and $\a \leq_F \be$;}
\item{(iii)}{$\a \in E, \ \be \in F$ and $\e(\a) \ \C \ \e'(\be)$.}
}
\enddefinition

\remark{Remark 1.1.6}
Definition 1.1.5 can easily be shown to be sound (see \cite{{\bf 7}, \S2}).
It is immediate from the construction that $E$ and $F$ are subheaps of 
$E \circ F$.

As in \cite{{\bf 7}}, we will write $\a \circ E$ and $E \circ \a$ for
$\{\a\} \circ E$ and $E \circ \{ \a \}$, respectively.  Note that $\a \circ E$
and $\be \circ E$ are equal as heaps if $\e(\a) = \e(\be)$.
\endremark

\definition{Definition 1.1.7}
The {\it concurrency graph} associated to the class of heaps $H(P, \C)$ is
the graph whose vertices are the elements of $P$ and for which there is an
edge from $v \in P$ to $w \in P$ if and only if $v \ne w$ and $v \ \C \ w$.
If $E = [E, \leq, \e]$ is a heap of $H(P, \C)$, we define the {\it concurrency
subgraph} of $E$ to be the full subgraph of the concurrency graph of $H(P, \C)$
that contains the vertices $\{\e(a) : a \in E\}$.
\enddefinition

\subhead 1.2 Rank functions \endsubhead

We now give our definition of the rank function and develop some of its 
elementary properties.

\definition{Definition 1.2.1}
Let $(E, \leq)$ be a poset.  If $a, b \in E$, the relation $a < b$ is
said to be a {\it covering relation} if there does not exist $c \in E$
such that $a < c < b$.
A function $\rho : E \ra \zed$ is said to be a 
{\it rank function} for $(E, \leq)$ if whenever $a, b \in E$ are such that
$a < b$ is a covering relation, we have $\rho(b) = \rho(a) + 1$.  If a
rank function for $(E, \leq)$ exists, we say $(E, \leq)$ is {\it ranked}.
\enddefinition

There are variants of Definition 1.2.1 in the literature, but our formulation 
is convenient for our purposes.

\definition{Definition 1.2.2}
Let $(E, \leq)$ be a poset and let $a, b \in E$.  We write $a \sim_c b$ if
$a < b$ is a covering relation, and we denote the equivalence relation on 
$E$ generated by $\sim_c$ by $\sim$.  We call the $\sim$-equivalence classes
of $E$ the {\it connected components} of $(E, \leq)$.
\enddefinition

The following lemma is clear from the definitions.

\proclaim{Lemma 1.2.3}
Let $(E, \leq)$ be a poset and let $\kappa : E \ra \zed$ be a function 
constant on $\sim$-equivalence classes.  If $\rho$ is a rank function for
$(E, \leq)$, then so is the function $\rho + \kappa$ defined by 
$(\rho + \kappa)(z) = \rho(z) + \kappa(z)$. \qed
\endproclaim

\definition{Definition 1.2.4}
Let $(E, \leq, \e)$ be a labelled heap.  We say $(E, \leq, \e)$ is ranked 
if the underlying poset $(E, \leq)$ is ranked.  In this case, we also say 
that the heap $[E, \leq, \e]$ is ranked.
\enddefinition

Definition 1.2.4 is sound because the property of being ranked is an invariant
of poset isomorphism.

\definition{Definition 1.2.5}
Let $(E, \leq)$ be a poset and let $a, b \in E$.  The 
interval $[a, b]$ is the subset
$\{x \in E : a \leq x \leq b\}$.  We make the same definition 
if $(E, \leq, \e)$ is a labelled heap.  If $[E, \leq, \e]$ is the 
corresponding heap, we call the subheap corresponding to the subset $[a, b]$
a {\it subinterval} of $[E, \leq, \e]$; we will often abuse notation and refer
to the subheap itself as $[a, b]$.  If $[a, b]$ is a subinterval in the heap 
$[E, \leq, \e]$, we say $[a, b]$ is a {\it balanced subinterval} if 
$\e(a) = \e(b)$.  A balanced subinterval $[a, b]$ is said to be {\it minimal}
if $a \ne b$ and if the only elements $c \in [a, b]$ with 
$\e(c) = \e(a) (= \e(b))$ are $c = a$ and $c = b$.
\enddefinition

We will regard subintervals of posets as subposets, in the obvious way.  The
following property will often be useful.

\remark{Remark 1.2.6}
If $[a, b]$ is a subinterval in $(E, \leq)$ and $x < y$ is a covering 
relation in the subinterval $[a, b]$ then $x < y$ is a covering relation 
in $(E, \leq)$.
\endremark

\proclaim{Lemma 1.2.7}
If $(E, \leq)$ is a ranked poset then every subinterval of $(E, \leq)$ is 
ranked.
\endproclaim

\demo{Proof}
Let $a, b \in E$ with $a < b$, and let $\rho$ be a rank function for 
$(E, \leq)$.  Then the restriction of $\rho$ to $[a, b]$ is a rank function
for $[a, b]$ by Remark 1.2.6.
\qed\enddemo

The main purpose of this paper is to investigate the extent to which the 
converse of Lemma 1.2.7 holds; we will see that the converse is false
in general.  The proof of the main result (Theorem 2.1.1) will
involve the following lemma.

\proclaim{Lemma 1.2.8}
Let $E = [E, \leq, \e]$ be a nonempty heap in $H(P, \C)$ and let 
$\a, \be \in E$.  If $\a \sim \be$ (as in Definition 1.2.2) then there 
is a sequence $$
\a = \g_0, \g_1, \ldots, \g_r = \be
$$ of elements of $E$ such that for each $0 \leq i < r$, we have
$\e(\g_i) \ \C \ \e(\g_{i+1})$.
\endproclaim

\demo{Proof}
Since $\a \sim \be$, the definition of $\sim$ shows that there is a (possibly
trivial) sequence
$\a = \g_0, \g_1, \ldots, \g_r = \be$ where, for each $0 \leq i < r$, 
either $\g_i < \g_{i+1}$ or $\g_i > \g_{i+1}$ is a covering relation.  
The lemma now follows from
part 2 of Definition 1.1.1.
\qed\enddemo

%\vfill\eject

\head 2. A sufficient condition for a heap to be ranked \endhead

We devote \S2 to investigating the converse of Lemma 1.2.7 for a general heap.
The main result of this section is Theorem 2.1.1.

\subhead 2.1. The main result \endsubhead

\proclaim{Theorem 2.1.1}
Let $E = [E, \leq, \e]$ be a heap in $H(P, \C)$.  Suppose the concurrency 
subgraph of $E$ (see Definition 1.1.7) contains no circuits.  Then the 
following are equivalent:
\item{\rm (i)}{$E$ is ranked;}
\item{\rm (ii)}{every subinterval of $E$ is ranked;}
\item{\rm (iii)}{every minimal balanced subinterval of $E$ is ranked.}
\endproclaim

\remark{Remark 2.1.2}
The implication (i) $\Rightarrow$ (ii) is immediate from Lemma 1.2.7 and 
the implication (ii) $\Rightarrow$ (iii) is trivial, so our strategy will
be to show that (iii) implies (i).
\endremark

\remark{Remark 2.1.3}
The circuit avoidance property above is called property (H4) in \cite{{\bf 6}}.
Some restriction is necessary here (see Example 2.1.5), although the 
condition given is too strong (see Example 2.1.4).
\endremark

\example{Example 2.1.4}
Let $P = \{1, 2, 3, 4, 5\}$ with concurrency relation $\C$ such that 
$a \ \C \ b$ for all $a, b \in P$; the concurrency graph $\Gamma$ 
is thus the complete graph on 5 vertices.  Let $E = [E, \leq, \e]$ be 
any of the heaps of $H(P, \C)$ with concurrency subgraph equal to $\Gamma$.
In this case, $(E, \leq)$ is totally ordered, and it follows that $E$ is
a ranked heap, as are all of its subintervals.
However, $\Gamma$ contains circuits.
\endexample

\example{Example 2.1.5}
Let $P = \{1, 2, 3, 4, 5\}$ as in Example 2.1.4, but define the concurrency
relation $\C$ so that $a \ \C \ b$ if and only if $\{a, b\}$ is in the list $$
\{ 
\{1\}, \{2\}, \{3\}, \{4\}, \{5\},
\{1, 2\}, \{2, 3\}, \{3, 4\},  \{4, 5\}, \{5, 1\}
\}.$$  In this case, $\Gamma$ is a pentagon.  Figure 1 shows the Hasse
diagram of a heap $E$ with concurrency subgraph $\Gamma$.  (This notation is
familiar from \cite{{\bf 6}}: for example we can see from the diagram that the
two minimal elements of $E$ are labelled $3$ and $5$, and the two 
maximal elements are labelled $1$ and $4$.)  It is not hard to see
that no rank 
function for $E$ exists, but that all subintervals of $E$ are ranked.
This is possible because the concurrency subgraph of $E$
contains a circuit.
\endexample

\topcaption{Figure 1} The heap $E$ of Example 2.1.5
\endcaption
\centerline{
\hbox to 1.069in{
\vbox to 0.902in{\vfill
        \includegraphics{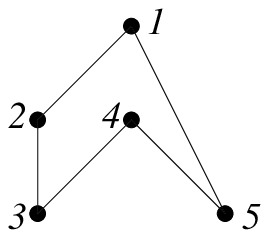}
}
\hfill}
}

\subhead 2.2 Proof of the main result \endsubhead

\proclaim{Lemma 2.2.1}
Let $E = [E, \leq, \e]$ be a nonempty heap in $H(P, \C)$, and let 
$\a$ be a minimal element of $E$.  Let $F$ be the subheap of $E$ 
corresponding to the subset $E\backslash \{\a\}$, so that $E = \a \circ F$.
Suppose that $F$ is ranked and that every minimal
balanced subinterval of $E$ is ranked, and suppose further that the concurrency
subgraph of $E$ contains no circuits.  If $\be, \g \in F$ are in the same
connected component of $F$ and $\a < \be$ and $\a < \g$ are covering 
relations in $E$, then we have $\rho(\be) = \rho(\g)$ for any rank function
$\rho$ of $F$.
\endproclaim

\demo{Proof}
We may assume that $F$ is not empty and that $\be \ne \g$, or there is 
nothing to prove.  Let $\Gamma$ be the concurrency subgraph of $E$; it contains
no circuits by hypothesis.
The condition $\be \ne \g$ and Definition 1.1.1 imply that the pieces
$\e(\be)$, $\e(\a)$ and $\e(\g)$ are distinct; since $\a < \be$ and 
$\a < \g$ are also covering relations, it must be the case that
$(\e(\be), \e(\a), \e(\g))$ is a sequence of distinct, adjacent vertices
in $\Gamma$.

By Lemma 1.2.8, there is a sequence $$
\be = \g_0, \g_1, \ldots, \g_r = \g
$$ of elements of $F = [F, \leq, \e]$ such that for each $0 \leq i <
r$, either $\e(\g_i) = \e(\g_{i+1})$ or $\e(\g_i)$ is adjacent to 
$\e(\g_{i+1})$ in $\Gamma$.  Since $\Gamma$
contains no circuits, the remarks in the first paragraph of the proof show
that every path from $\e(\be)$ to $\e(\g)$ passes through $\e(\a)$, and 
therefore $\e(\g_i) = \e(\a)$ for some $0 < i < r$.  This means that 
there is an element $\a' \in F$ with $\e(\a') = \e(\a)$.  

The subinterval $[\a, \a']$ of $E$ is balanced, and
so $E$ contains a minimal balanced subinterval $[\a, \a'']$ for some 
$\a'' \in F$.  Now $\a''$ is 
comparable to both $\be$ and $\g$ in the partial order, and condition 1
of Definition 1.1.1 implies that $\be < \a''$ and $\g < \a''$.
Since $\be \in [\a, \a'']$, there must be a sequence $$
\be = \be_0 < \be_1 < \cdots < \be_t = \a''
$$ where each of the relations $\be_i < \be_{i+1}$ is a covering relation in
$[\a, \a'']$, and therefore (by Remark 1.2.6) in $E$.  

Note that $[\a, \a'']$ is ranked as a subinterval of $E$ by
hypothesis; this implies that the saturated chains from $\a$ to $\a''$
have a common length.
Fixing a rank function $\rho$ for $F$, we now find that
$\rho(\a'') = \rho(\be) + t$; similarly, 
$\rho(\a'') = \rho(\g) + t'$, where $t'$ is the length of a saturated
chain from $\g$ to $\a''$.  (Note that $t$ and $t'$ are independent
of $\rho$.)  Because $\a < \be$ and $\a < \g$ are covering relations,
the above assertion about saturated chains forces $t = t'$, and we
have $\rho(\be) = \rho(\g)$ as required.
\qed\enddemo

\demo{Proof of Theorem 2.1.1}
By Remark 2.1.2, it is enough to prove the implication (iii) $\Rightarrow$ (i).
Let $E = [E, \leq, \e]$ be a heap in $H(P, \C)$.  Suppose the concurrency 
subgraph of $E$ contains no circuits and that every minimal balanced 
subinterval of $E$ is ranked.  The proof is by induction on $|E|$.  If 
$|E|$ is $0$ or $1$, $E$ will be ranked for trivial reasons and there is
nothing to prove.  We may therefore assume that $E = \a \circ F$ for some
subheap $F$ of $E$ with $|F| = |E| - 1$, and suppose that $\rho$ is a rank
function for $F$.  (It is clear that all subheaps of $E$ will also have 
concurrency graphs with no circuits.)

If $\a$ is the only element in its connected component in $E$, we may extend
$\rho$ to $E$ by defining $\rho(\a)$ arbitrarily.  Otherwise, since $\a$ is
minimal in $E$, we have covering relations $\a < \be_i$ for some nonempty
set $\{\be_i\} \subset F$.  If $\be_i$ and $\be_j$ are in the same connected
component of $F$ then Lemma 2.2.1 shows that $\rho(\be_i) = \rho(\be_j)$ 
By using Lemma 1.2.3 (if necessary) to adjust the values of the rank 
function on the connected components of $F$, we may assume that $\rho$
is constant on the set $\{\be_i\}$.  The proof is completed by defining
$\rho(\a) := \rho(\be) - 1$ for (any) $\be \in \{\be_i\}$.  
\qed\enddemo

\head 3. Heaps of fully commutative elements in Coxeter groups \endhead

In \S3, we turn our attention to the special case of heaps that arise
from fully commutative elements of Coxeter groups; these were studied
by Stembridge in \cite{{\bf 5}}.  It turns out (Theorem 3.2.3) 
that if we restrict our attention
to Coxeter groups having only finitely many fully commutative elements,
it becomes easy to determine whether every minimal balanced
subinterval of the heap is ranked.  The result does not hold if we drop
the finiteness hypothesis, and the proof relies on the classification of
such Coxeter groups, but it is nevertheless potentially very 
helpful when checking examples by hand or by computer.

\vfill\eject

\subhead 3.1 Heaps of fully commutative elements \endsubhead

\definition{Definition 3.1.1}
A {\it Coxeter group} is a pair $(W, S)$ where $S$ is a set and $W$ is the
group generated by $S$ subject to the defining relations $$
(st)^{m(s, t)} = 1
,$$ where $m(s, s) = 1$ for $s \in S$ and 
$2 \leq m(s, t) = m(t, s) \leq \infty$ for $s, t \in S$ and $s \ne t$.
(For the purposes of this paper, we will always assume that the set 
$S$ is finite.)
The {\it Coxeter graph} of $(W, S)$ has vertex set $S$.  Two 
distinct vertices $s, t$ in the Coxeter graph are joined by an edge labelled
$m = m(s, t)$ if $m \geq 3$, but if $m = 3$ we omit the label on the
edge by convention.
\enddefinition

We take the following to be the definition of the heap of a fully commutative
element; this is not the original definition but is equivalent to it by
\cite{{\bf 5}, Proposition 2.3}.  In this paper, we are not concerned with the 
fully commutative elements of Coxeter groups themselves, but rather only
with their heaps.  

\definition{Definition 3.1.2}
Let $(W, S)$ be a Coxeter group.  We define $\C$ by the condition $$
s \ \C \  t \Leftrightarrow m(s, t) \ne 2
.$$  A heap $E = [E, \leq, \e]$ in $H(S, \C)$ is the {\it heap of a fully 
commutative element of $W$} if and only if the following conditions hold.

\item{1.}{There is no convex chain $\a_1 < \a_2 < \cdots < \a_m$ in $E$
such that $\e(\a_i) = s$ for all odd $i$ and $\e(\a_i) = t$ for all even $i$,
where $3 \leq m = m(s, t) < \infty$.}
\item{2.}{There is no covering relation $\a < \be$ in $E$ such that
$\e(\a) = \e(\be)$.}

We say $(W, S)$ is an {\it FC-finite Coxeter group} if the number of (heaps
of) fully commutative elements is finite.
\enddefinition

\remark{Remark 3.1.3}
The fully commutative elements of $W$ are in bijection with heaps satisfying
the conditions of Definition 3.1.2; for an explanation see \cite{{\bf 5}, \S1.2}.
\endremark

\remark{Remark 3.1.4}
The term ``convex chain'' in Definition 3.1.2 has its obvious meaning: a
chain $$
\be_1 < \be_2 < \cdots < \be_r
$$ in $E$ is said to be convex if, whenever $\g \in E$ is such that 
$\be_i < \g < \be_j$ for some $1 \leq i, j \leq r$, $\g$ lies in the chain.
\endremark

\example{Example 3.1.5}
Consider a Coxeter graph of type $D_5$, meaning that $$
S = \{1, 2, 3, 4, 5\}$$ and $m(s, t) = 2$ 
unless $s = t$ (in which case $m(s, t) = 1$) 
or $\{s, t\}$ is one of the pairs $$
\{1, 3\}, \{2, 3\}, \{3, 4\}, \{4, 5\}
$$ (in which case $m(s, t) = 3$).  Figure 2 shows a fully commutative heap of
type $D_5$, \idest of type $H(S, \C)$ where $\C$ is as in Definition 3.1.2.
The (unique) chain corresponding to the sequence of labels $(3, 1, 3)$ is
not convex, due to the position of the occurrence of the label $2$.  One checks
similarly that there are no chains violating condition 1 of Definition 3.1.2.
It is easy to verify that the situation in condition 2 of 
Definition 3.1.2 cannot occur.
\endexample

\topcaption{Figure 2} A fully commutative heap of type $D_5$
\endcaption
\centerline{
\hbox to 1.027in{
\vbox to 1.680in{\vfill
        \includegraphics{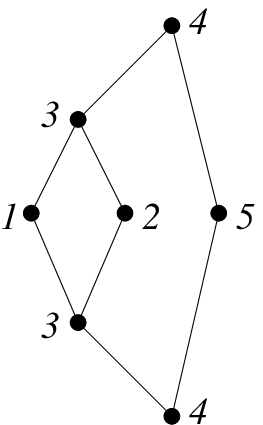}
}
\hfill}
}

The classification of FC-finite Coxeter groups in terms of their Coxeter
graphs was given by Stembridge \cite{{\bf 5}, Theorem 4.1}, and a similar 
result was independently obtained by Graham \cite{{\bf 2}, Theorem 7.1} 
from an algebraic perspective.

\proclaim{Theorem 3.1.6 (Stembridge; Graham)}
A Coxeter group $(W, S)$ is FC-finite if and only if the connected components
of its Coxeter graph appear in the list in Figure 3. \qed
\endproclaim

%\vfill\eject

\topcaption{Figure 3} Connected components of Coxeter graphs of 
FC-finite Coxeter groups
\endcaption
\centerline{
\hbox to 3.319in{
\vbox to 4.194in{\vfill
        \includegraphics{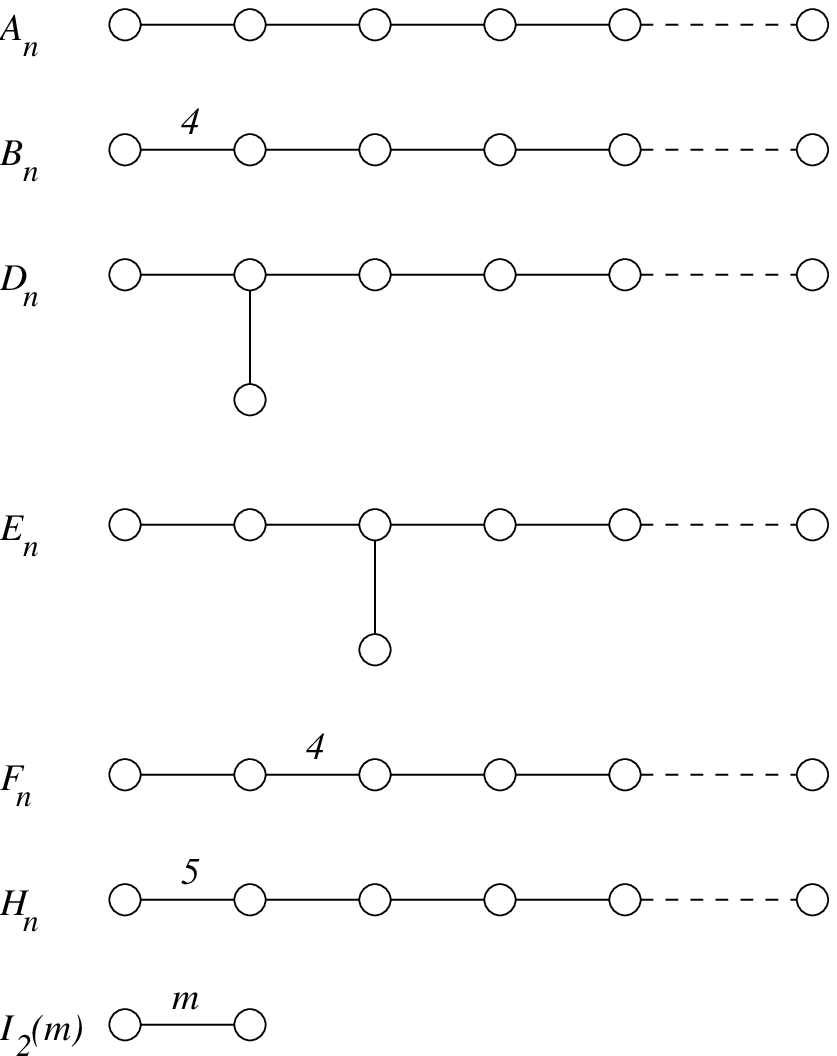}
}
\hfill}
}

(The subscript $n$ in Figure 3 denotes the number of vertices in the graph.)

\subhead 3.2 Ranked heaps of fully commutative elements \endsubhead

The main result of \S3 is Theorem 3.2.3, which gives a concise 
characterization of ranked heaps of fully commutative elements in FC-finite
Coxeter groups.

\definition{Definition 3.2.1}
Let $E = [E, \leq, \e]$ be a heap in $H(P, \C)$, and let $[a, b]$ be a minimal
balanced subinterval of $E$.  We define the subset $S_{[a, b]}$ of $E$ by $$
S_{[a, b]} = \{ c \in [a, b] : \e(a) \ne \e(c) \text{ and } \e(a) \ \C \ 
\e(c)\}.
$$
\enddefinition

\example{Example 3.2.2}
Let $a$ and $b$ be the minimal and maximal elements of the heap shown in 
Figure 2.  Then $S_{[a, b]}$ consists of three elements of the subinterval
$[a, b]$ (which is in this case the whole heap): the one
labelled $5$ and the two labelled $3$.
\endexample

\proclaim{Theorem 3.2.3}
Let $(W, S)$ be an FC-finite Coxeter group, and let $E$ be the heap of
a fixed fully commutative element $w \in W$.  The following are equivalent:
\item{\rm (i)}{$E$ is ranked;}
\item{\rm (ii)}{for each minimal balanced subinterval $[a, b]$ of $E$, either
(a) all the elements of $S_{[a, b]}$ have the same label or
(b) all the elements of $S_{[a, b]}$ have distinct labels.}
\endproclaim

\example{Example 3.2.4}
Let $a$ and $b$ be the minimal and maximal elements of the heap in Figure 2.
Theorem 3.2.3 applies because a Coxeter group of type $D_5$ is FC-finite by
Theorem 3.1.6, and the heap in question corresponds to a fully commutative 
element by Example 3.1.5.
The three elements of $S_{[a, b]}$ do not all have the same label, but they do
not have distinct labels either, so the heap is not ranked.
\endexample

\subhead 3.3 Proof of Theorem 3.2.3 \endsubhead

\definition{Definition 3.3.1}
Let $\Gamma$ be a Coxeter graph and let $s$ and $t$ be adjacent vertices of
$\Gamma$.  Let $\Gamma\backslash\{s\}$ be the graph obtained from $\Gamma$
by deleting $s$ and all edges emerging from $s$, let $\Gamma_{s, t}$ be the 
connected component of $\Gamma\backslash\{s\}$ that contains $t$, and let
$\Gamma_{s \sra t}$ be the full subgraph of $\Gamma$ containing $s$ and
the vertices of $\Gamma_{s, t}$.
\enddefinition

\example{Example 3.3.2}
Let $\Gamma$ be a graph of type $E_8$ as shown in Figure 3, let $s$ be the
vertex of degree $3$ and let $t$ be the vertex immediately to the right of
$s$.  Then $\Gamma\backslash\{s\}$ consists of the disjoint union of three
Coxeter graphs of types $A_1$, $A_2$ and $A_4$; $\Gamma_{s, t}$ is a Coxeter
graph of type $A_4$ and $\Gamma_{s \sra t}$ is a Coxeter graph of type $A_5$
containing $s$ and all the vertices to the right of $s$.
\endexample

\proclaim{Lemma 3.3.3}
Let $\Gamma$ be the Coxeter graph of an FC-finite Coxeter group and let $s$
be a vertex of $\Gamma$ with degree strictly greater than $1$.
There is at most one vertex $t$ adjacent to $s$ such that 
$\Gamma_{s \sra t}$ is not of type $A_n$ for some $n \geq 2$.
\endproclaim

\demo{Proof}
This is a case by case check using Theorem 3.1.6 (see Figure 3).
\qed\enddemo

\proclaim{Lemma 3.3.4}
Let $E = [E, \leq, \e]$ be the heap (in $H(P, \C)$) 
of a fully commutative element in an FC-finite Coxeter group
and let $[a, b]$ be a minimal balanced subinterval of $E$.  Suppose the 
elements of $S_{[a, b]}$ do not all have the same label.  Then there exists
an element of $S_{[a, b]}$ whose label is unique among the labels of elements
of $S_{[a, b]}$.
\endproclaim

\demo{Proof}
Since the elements of $S_{[a, b]}$ do not all have the same label, the
degree of $\e(a)$ in the concurrency graph $\Gamma$ is greater than 1.  
Let $c, d \in S_{[a, b]}$ be such that $\e(c) \ne \e(d)$; both labels are 
distinct from $\e(a) = \e(b)$ by minimality of the subinterval.  By
Lemma 3.3.3, we may assume
without loss of generality that $\Gamma_{\e(a) \sra \e(c)}$ is of type
$A_n$ for some $n \geq 2$.  We index the vertices of this subgraph of type $A$ 
by $p_1 = \e(a)$, $p_2 = \e(c)$, $p_3, \ldots, p_n$ such that $p_i$ and $p_j$
are adjacent in $\Gamma$ if and only if $|i - j| = 1$.

Suppose, for a contradiction, that $E$ is the heap of a fully
commutative element in an FC-finite Coxeter group, and that $[a, b]$ is a
minimal balanced subinterval of $E$ for which (a) the elements of $S_{[a,
b]}$ do not all have the same label and (b) there is no element of
$S_{[a, b]}$ whose label is unique among the labels of elements in
$S_{[a, b]}$.
We claim by induction that for each $1 \leq k < n$, there is a minimal 
balanced subinterval $[a_k, b_k]$ with $\e(a_k) = \e(b_k) = p_k$ 
containing at least two elements labelled
$p_{k+1}$.  Define $a = a_1$, $b = b_1$ and 
observe that $[a, b]$ contains at least one element labelled $p_2$ by 
definition of $c$.  By part (b) of the assertion above, there must be
at least two elements of $[a, b]$ labelled $p_2$, which
establishes the $k = 1$ case of the induction.

For the inductive step, we may assume $n > 2$.  Suppose $k < n - 1$ and 
that $[a_k, b_k]$ contains 
at least two elements, $a'$ and $b'$, labelled $p_{k+1}$.  We may assume
that the balanced chain $[a', b']$ is minimal by choosing $a'$ and $b'$
suitably.  By minimality of $[a_k, b_k]$, we see that $[a', b']$ contains
no elements labelled $p_k$.  Since $\Gamma_{\e(a) \sra \e(c)}$ is of type
$A_n$, we must have at least two elements in $[a', b']$ labelled
$p_{k+2}$: if there were none, we would have a counterexample to 
condition 2 of Definition 3.1.2 by taking $\a = a', \be = b'$, and if there
were only one, we would have a counterexample to condition 1 of that
definition by taking $\a_1 = a', \a_3 = b'$ and $\a_2$ to be the element
labelled $p_{k+2}$.  This proves the inductive step after taking
$a_{k+1} = a'$, $b_{k+1} = b'$.

This situation leads to a contradiction because $[a_n, b_n]$ is a minimal
balanced subinterval containing no occurrences of $p_{n-1}$ 
(using the case $k = n - 1$ above).  Taking $\a = a_n, \be = b_n$
in condition 2 of Definition 3.1.2 shows that $E$ is not the heap of a fully
commutative element, a contradiction.
\qed\enddemo

\proclaim{Lemma 3.3.5}
Let $E = [E, \leq, \e]$ be a heap in $H(P, \C)$ such that the concurrency
subgraph of $E$ contains no circuits, and let $[a, b]$ be a minimal
balanced subinterval of $E$.  Suppose $c \in S_{[a, b]}$ and define $a'$ 
(respectively, $b'$) to be the minimal (respectively, maximal) element of
$S_{[a, b]}$ with label $\e(c)$.  Then $a < a'$ and $b' < b$ are covering
relations in $E$.
\endproclaim

\demo{Proof}
We deal with the case of $a'$; the other case is similar.  Since
$\e(a) \ \C \ \e(a')$, there is a chain of covering relations $$
a = a_0 < a_1 < \cdots < a_t = a'
.$$  The definition of $a'$ ensures that $t > 0$, and we are done if
$t = 1$, so suppose $t > 1$.
Since $a' < b$, minimality of $[a, b]$ shows that if $i > 0$ then
$a_i$ cannot have label $\e(a)$.  Similarly, the definition of $a'$
shows that if $i < t$ then $a_i$ cannot have label $\e(a')$.
By Lemma 1.2.8, the corresponding sequence $$
\e(a_0), \e(a_1), \ldots, \e(a_t)
$$ in $P$ is a path (possibly with repeated vertices) between the adjacent
vertices $\e(a_0)$ and $\e(a_t)$ that passes through each of $\e(a_0)$ and
$\e(a_t)$ precisely once, which is impossible as $t > 1$ 
and the concurrency graph contains no circuits.  This completes the proof.
\qed\enddemo

\example{Example 3.3.6}
Maintain the set-up in Example 3.2.2; recall that this concerns 
the heap of a fully commutative element.  As noted in Example 3.2.2,
the elements
of $S_{[a, b]}$ do not all have the same label; Lemma 3.3.4 then predicts
that one of the labels ($5$ in this case) occurs uniquely in the subinterval
$[a, b]$.  (This is because $\Gamma_{4 \sra 5}$ is of type $A_2$.)  Lemma
3.3.5 predicts that each of the elements labelled $3$ or $5$  covers
or is covered by either $a$ or $b$.
\endexample

\demo{Proof of Theorem 3.2.3}
Since $(W, S)$ is an FC-finite Coxeter group, the concurrency graph of $E$ has
no circuits because none of the graphs in Figure 3 has any circuits.  (The
relation between the Coxeter graph and the concurrency graph is given in
Definition 3.1.2.)  

First, suppose $E$ is ranked.  By Theorem 2.1.1, every minimal
balanced subinterval of $E$ is ranked; 
let $[a, b]$ be such an subinterval.  If all the elements of $S_{[a, b]}$ have
the same label then condition (ii) of Theorem 3.2.3 holds, and we are
done.  If not, Lemma 3.3.4 shows the existence of an element $c \in [a, b]$ 
whose label is unique among the labels of $S_{[a, b]}$.  By Lemma 3.3.5,
$a < c$ and $c < b$ are covering relations, which means that if $\rho$ is
any rank function for $E$ then $\rho(b) = \rho(a) + 2$.  Suppose the statement
of Theorem 3.2.3 (ii) does not hold, so that there exist at least two elements
$d, d' \in S_{[a, b]}$ with $\e(d) = \e(d') \ne \e(c)$.  Without loss of
generality, $d < d'$, so we have a chain $a < d < d' < b$.  This means that
$\rho(b) > \rho(a) + 2$, a contradiction, and condition (ii) of Theorem 3.2.3
holds, as required.

For the converse, we will prove by induction on $|E|$ that (ii) implies (i).
If $|E|$ is $0$ or $1$ the heap $E$ is ranked for trivial reasons and there
is nothing to prove.  For the general case, assume the hypotheses of (ii) and
consider an arbitrary minimal balanced subinterval $[a, b]$ in $E$.  
If we can prove that $[a, b]$ is ranked, the claim will follow by 
Theorem 2.1.1.  There are two cases to consider.

In the first case, the labels of the elements
$c_1, c_2, \ldots, c_r$ of $S_{[a, b]}$ are distinct.  Lemma 3.3.5 shows
that $a < c_i < b$ is a chain of covering relations for each $i$, so the
subinterval $[a, b]$ consists only of the elements $c_i$ together with $a$ and
$b$.  The subinterval is ranked in this case: we may take $\rho(a) = 0$, 
$\rho(b) = 2$ and $\rho(c_i) = 1$ for each $i$.

In the second case to be considered, the elements $c_1, c_2, \ldots, 
c_r$ of $S_{[a, b]}$ all have the same label, so we may assume that
$c_1 < c_2 < \cdots < c_r$.  By Lemma 3.3.5, $a < c_1$ and $c_r < b$ are
covering relations in $E$; there are no other covering relations of
the form $a < c'$ or $c' < b$ by the assumption on $S_{[a, b]}$.  
It follows that the subinterval $[a, b]$ consists (as a set) of the 
balanced subinterval $[c_1, c_r]$ together with the additional 
elements $a$ and $b$.  We claim that any subinterval in the heap of a fully
commutative element is also the heap of a fully commutative element for
the same Coxeter group: this follows from Definition 3.1.2 and the general 
fact that any convex chain in a subinterval of a poset is also a convex chain 
in the poset.  Furthermore, we claim that any minimal balanced subinterval 
$[d, e]$ of an subinterval in a heap $E$ is also a minimal balanced 
subinterval of $E$: it is minimal because the set of elements in $E$
with a given label is totally ordered.  These two observations show
that $[c_1, c_r]$ is the heap of a fully commutative element $w \in
W$, and that it satisfies condition (ii) of
Theorem 3.2.3.  The subinterval $[c_1, c_r]$ contains strictly fewer elements
than $E$ and is therefore ranked by the inductive hypothesis; let $\rho$ be
a rank function for $[c_1, c_r]$.  We can extend $\rho$ to a rank function for
$[a, b]$ by defining $\rho(a) = \rho(c_1) - 1$ and $\rho(b) = \rho(c_r) + 1$.
\qed\enddemo

\remark{Remark 3.3.7}
If $E$ is a heap of fully commutative element of a Coxeter group of
type $A_n$, it is well known and easy to show using the techniques of the
proof of Lemma 3.3.4 that if $[a, b]$ is a minimal balanced subinterval of
$E$ then $S_{[a, b]}$ consists of precisely two elements, with distinct labels.
It follows that any heap of a fully commutative element of a Coxeter group
of type $A_n$ is ranked.  This is also well known and is what allows Billey 
and Warrington's method of ``pushing together the connected components of
a heap'' \cite{{\bf 1}, \S3} to work.
\endremark

\head Acknowledgements \endhead

The author thanks J.R. Stembridge for helpful correspondence, and the
referee for suggesting many improvements to an earlier version of this paper.

\leftheadtext{}
\rightheadtext{}
%\vfill\eject

\end
\Refs\refstyle{A}\widestnumber\key{{\bf 7}}
\leftheadtext{References}
\rightheadtext{References}

\ref\key{{\bf 1}}
\by S.C. Billey and G.S. Warrington
\paper Kazhdan--Lusztig Polynomials for 321-hexagon-avoiding permutations
\jour J. Algebraic Combin.
\vol 13 \yr 2001 \pages 111--136
\endref

\ref\key{{\bf 2}}
\by J.J. Graham
\book Modular representations of Hecke algebras and related algebras
\publ Ph.D. thesis
\publaddr University of Sydney
\yr 1995
\endref

\ref\key{{\bf 3}}
\by A. Mazurkiewicz
\paper Trace theory
\inbook Petri nets, applications and relationship to other models of
concurrency
\publ Lecture Notes in Computer Science 255, Springer
\publaddr Berlin--Heidelberg--New York
\yr 1987
\pages 279--324
\endref

\ref\key{{\bf 4}}
\by R.A. Proctor
\paper Minuscule elements of Weyl groups, the numbers game, and
$d$-complete posets
\jour J. Algebra
\vol 213 \yr 1999 \pages 272--303
\endref

\ref\key{{\bf 5}}
\by J.R. Stembridge 
\paper On the fully commutative elements of Coxeter groups 
\jour J. Algebraic Combin.
\vol 5 
\yr 1996 
\pages 353--385
\endref

\ref\key{{\bf 6}}
\by J.R. Stembridge
\paper Minuscule elements of Weyl groups
\jour J. Algebra 
\vol 235 \yr 2001 \pages 722--743
\endref

\ref\key{{\bf 7}}
\by G.X. Viennot
\paper Heaps of pieces, I: basic definitions and combinatorial lemmas
\inbook Combinatoire \'E\-nu\-m\'e\-ra\-tive
\publ Springer-Verlag
\publaddr Berlin
\yr 1986 \pages 321--350 \bookinfo ed. G. Labelle and P. Leroux
\endref

\endRefs
